\documentclass{amsart}
\usepackage{amsmath}
\usepackage{amsfonts}

\theoremstyle{plain}

\numberwithin{equation}{section}
\input{tcilatex}

\begin{document}
\title[ON WEIGHTED $q$-GENOCCHI NUMBERS AND POLYNOMIALS]{\textbf{ANALYTIC
CONTINUATION OF WEIGHTED }$q$\textbf{-GENOCCHI NUMBERS AND POLYNOMIALS}}
\author[\textbf{S. Arac\i }]{\textbf{Serkan Arac\i }}
\address{University of Gaziantep, Faculty of Science and Arts, Department of
Mathematics, 27310 Gaziantep, TURKEY}
\email{mtsrkn@hotmail.com}
\author[\textbf{M. Acikgoz}]{\textbf{Mehmet Acikgoz}}
\address{University of Gaziantep, Faculty of Science and Arts, Department of
Mathematics, 27310 Gaziantep, TURKEY}
\email{acikgoz@gantep.edu.tr\bigskip \bigskip }
\author[\textbf{A. G\"{u}rsul}]{\textbf{Aynur G\"{u}rsul}}
\address{University of Gaziantep, Faculty of Science and Arts, Department of
Mathematics, 27310 Gaziantep, TURKEY}
\email{aynurgursul@hotmail.com}

\begin{abstract}
In the present paper, we analyse analytic continuation of weighted $q$%
-Genocchi numbers and polynomials. A novel formula for weighted $q$%
-Genocchi-Zeta function $\widetilde{\zeta }_{G,q}\left( s\mid \alpha \right) 
$ in terms of nested series of $\widetilde{\zeta }_{G,q}\left( n\mid \alpha
\right) $ is derived. Moreover, we introduce a novel concept of dynamics of
the zeros of analytically continued weighted $q$-Genocchi polynomials.
\end{abstract}

\maketitle

\section{\textbf{INTRODUCTION}}

In this paper, we use notations like $%
\mathbb{N}
$, $%
\mathbb{R}
$ and $%
\mathbb{C}
$, where $%
\mathbb{N}
$ denotes the set of natural numbers, $%
\mathbb{R}
$ denotes the field of real numbers and $%
\mathbb{C}
$ also denotes the set of complex numbers. When one talks of $q$-extension, $%
q$ is variously considered as an indeterminate, a complex number or a p-adic
number.

Throughout this work, we will assume that $q\in 
\mathbb{C}
$ with $\left\vert q\right\vert <1$. The $q$-integer symbol $\left[ x:q%
\right] $ denotes as 
\begin{equation*}
\left[ x:q\right] =\frac{q^{x}-1}{q-1}\text{.}
\end{equation*}

Firstly, analytic continuation of $q$-Euler numbers and polynomials was
investigated by Kim in \cite{Kim 1}. He gave a new concept of dynamics of
the zeros of analytically continued $q$-Euler polynomials. Actually, we were
motivated from his excellent paper which is "Analytic continuation of $q$%
-Euler numbers and polynomials, Applied Mathematics Letters 21 (2008)
1320-1323." We also procure to \ analytic continuation of weighted $q$%
-Genocchi numbers and polynomials as parallel to his article. Also, we give
some interesting identities by using generating function of weighted $q$%
-Genocchi polynomials.

\section{\textbf{PROPERTIES OF THE WEIGHTED }$q$\textbf{-GENOCCHI NUMBERS
AND POLYNOMIALS}}

For $\alpha \in 
\mathbb{N}
\tbigcup \left\{ 0\right\} $, the weighted $q$-Genocchi polynomials are
defined by means of the following generating function:

For $x\in 
\mathbb{C}
$, 
\begin{equation}
\sum_{n=0}^{\infty }\widetilde{G}_{n,q}\left( x\mid \alpha \right) \frac{%
t^{n}}{n!}=\left[ 2:q\right] t\sum_{n=0}^{\infty }\left( -1\right)
^{n}q^{n}e^{t\left[ n+x:q^{\alpha }\right] }\text{.}  \label{Equation 1}
\end{equation}

As a special case $x=0$ into (\ref{Equation 1}), $\widetilde{G}_{n,q}\left(
0\mid \alpha \right) :=\widetilde{G}_{n,q}\left( \alpha \right) $ are called
weighted $q$-Genocchi numbers. By (\ref{Equation 1}), we readily derive the
following 
\begin{equation}
\frac{\widetilde{G}_{n+1,q}\left( x\mid \alpha \right) }{n+1}=\frac{\left[
2:q\right] }{\left[ \alpha :q\right] ^{n}\left( 1-q\right) ^{n}}%
\sum_{l=0}^{n}\binom{n}{l}\left( -1\right) ^{l}\frac{q^{\alpha lx}}{%
1+q^{\alpha l+1}}\text{,}  \label{Equation 2}
\end{equation}

where $\binom{n}{l}$ is the binomial coefficient. By expression (\ref%
{Equation 1}), we see that%
\begin{equation}
\widetilde{G}_{n,q}\left( x\mid \alpha \right) =q^{-\alpha x}\left(
q^{\alpha x}\widetilde{G}_{q}\left( \alpha \right) +\left[ x:q^{\alpha }%
\right] \right) ^{n}\text{,}  \label{Equation 3}
\end{equation}

with the usual convention of replacing $\left( \widetilde{G}_{q}\left(
\alpha \right) \right) ^{n}$ by $\widetilde{G}_{n,q}\left( \alpha \right) $
is used (for details, see \cite{Araci 1}, \cite{Araci 2}).

Let $\widetilde{T}_{q}^{\left( \alpha \right) }\left( x,t\right) $ be the
generating function of weighted $q$-Genocchi polynomials as follows:%
\begin{equation}
\widetilde{T}_{q}^{\left( \alpha \right) }\left( x,t\right)
=\sum_{n=0}^{\infty }\widetilde{G}_{n,q}\left( x\mid \alpha \right) \frac{%
t^{n}}{n!}\text{.}  \label{Equation 4}
\end{equation}

Then, we easily notice that%
\begin{equation}
\widetilde{T}_{q}^{\left( \alpha \right) }\left( x,t\right) =\left[ 2:q%
\right] t\sum_{n=0}^{\infty }\left( -1\right) ^{n}q^{n}e^{t\left[
n+x:q^{\alpha }\right] }\text{.}  \label{Equation 5}
\end{equation}

From expressions (\ref{Equation 4}) and (\ref{Equation 5}), we procure the
followings:

For $k$ (=even) and $n,\alpha \in 
\mathbb{N}
\tbigcup \left\{ 0\right\} $, we have 
\begin{equation}
q^{k}\frac{\widetilde{G}_{n+1,q}\left( k\mid \alpha \right) }{n+1}-\frac{%
\widetilde{G}_{n+1,q}\left( \alpha \right) }{n+1}=\left[ 2:q\right]
\sum_{l=0}^{k-1}\left( -1\right) ^{l}q^{k-l-1}\left[ l:q^{\alpha }\right]
^{n}\text{.}  \label{Equation 6}
\end{equation}

For $k$ (=odd) and $n,\alpha \in 
\mathbb{N}
\tbigcup \left\{ 0\right\} $, we have 
\begin{equation}
q^{k}\frac{\widetilde{G}_{n+1,q}\left( k\mid \alpha \right) }{n+1}+\frac{%
\widetilde{G}_{n+1,q}\left( \alpha \right) }{n+1}=\left[ 2:q\right]
\sum_{l=0}^{k-1}\left( -1\right) ^{l}q^{k-l-1}\left[ l:q^{\alpha }\right]
^{n}\text{.}  \label{Equation 7}
\end{equation}

Via Eq. (\ref{Equation 5}), we easily obtain the following:%
\begin{equation}
\widetilde{G}_{n,q}\left( x\mid \alpha \right) =q^{-\alpha x}\sum_{k=0}^{n}%
\binom{n}{k}q^{\alpha kx}\widetilde{G}_{k,q}\left( \alpha \right) \left[
x:q^{\alpha }\right] ^{n-k}\text{.}  \label{Equation 8}
\end{equation}

From (\ref{Equation 6})-(\ref{Equation 8}), we get the following:%
\begin{eqnarray}
&&\left[ 2:q\right] \sum_{l=0}^{k-1}\left( -1\right) ^{l}q^{k-l-1}\left[
l:q^{\alpha }\right] ^{n}  \label{Equation 9} \\
&=&\left( q^{\alpha kn}-1\right) \frac{\widetilde{G}_{n+1,q}\left( \alpha
\right) }{n+1}+q^{-\alpha k}\sum_{j=0}^{n}\frac{1}{n+1}\binom{n+1}{j}%
q^{\alpha jk}\widetilde{G}_{k,q}\left( \alpha \right) \left[ k:q^{\alpha }%
\right] ^{n+1-k}\text{,}  \notag
\end{eqnarray}

here $k$ is an even positive integer. If $k$ is an odd positive integer.
Then, we can derive the following equality:%
\begin{eqnarray}
&&\left[ 2:q\right] \sum_{l=0}^{k-1}\left( -1\right) ^{l}q^{k-l-1}\left[
l:q^{\alpha }\right] ^{n}  \label{Equation 10} \\
&=&\left( q^{\alpha kn}+1\right) \frac{\widetilde{G}_{n+1,q}\left( \alpha
\right) }{n+1}+q^{-\alpha k}\sum_{j=0}^{n}\frac{1}{n+1}\binom{n+1}{j}%
q^{\alpha jk}\widetilde{G}_{k,q}\left( \alpha \right) \left[ k:q^{\alpha }%
\right] ^{n+1-k}\text{.}  \notag
\end{eqnarray}

\section{\textbf{WEIGHTED }$q$\textbf{-GENOCCHI-ZETA FUNCTION}}

The famous Genocchi polynomials were defined as%
\begin{equation}
\frac{2t}{e^{t}+1}e^{xt}=\sum_{n=0}^{\infty }G_{n}\left( x\right) \frac{t^{n}%
}{n!},\text{ }\left\vert t\right\vert <\pi \text{ cf. \cite{kim 4}.}
\label{Equation 11}
\end{equation}

For $s\in 
\mathbb{C}
$, $x\in 
\mathbb{R}
$ with $0\leq x<1$, Genocchi-Zeta function are given by%
\begin{equation}
\zeta _{G}\left( s,x\right) =2\sum_{n=0}^{\infty }\frac{\left( -1\right) ^{n}%
}{\left( n+x\right) ^{s}}\text{, }  \label{Equation 12}
\end{equation}

and%
\begin{equation}
\zeta _{G}\left( s\right) =\sum_{n=1}^{\infty }\frac{\left( -1\right) ^{n}}{%
n^{s}}\text{.}  \label{Equation 13}
\end{equation}

By (\ref{Equation 11}), (\ref{Equation 12}) and (\ref{Equation 13}),
Genocchi-Zeta functions are related to the Genocchi numbers as follows: 
\begin{equation*}
\zeta _{G}\left( -n\right) =\frac{G_{n+1}}{n+1}\text{.}
\end{equation*}

Moreover, it is simple to see%
\begin{equation*}
\zeta _{G}\left( -n,x\right) =\frac{G_{n+1}\left( x\right) }{n+1}\text{.}
\end{equation*}

The weighted $q$-Genocchi Hurwitz-Zeta type function are defined by 
\begin{equation*}
\widetilde{\zeta }_{G,q}\left( s,x\mid \alpha \right) =\left[ 2:q\right]
\sum_{m=0}^{\infty }\frac{\left( -1\right) ^{m}q^{m}}{\left[ m+x:q^{\alpha }%
\right] ^{s}}\text{ .}
\end{equation*}

Similarly, weighted $q$-Genocchi-Zeta function are given by%
\begin{equation*}
\widetilde{\zeta }_{G,q}\left( s\mid \alpha \right) =\left[ 2:q\right]
\sum_{m=1}^{\infty }\frac{\left( -1\right) ^{m}q^{m}}{\left[ m:q^{\alpha }%
\right] ^{s}}\text{.}
\end{equation*}

For $n,\alpha \in 
\mathbb{N}
\tbigcup \left\{ 0\right\} $, we have%
\begin{equation*}
\widetilde{\zeta }_{G,q}\left( -n\mid \alpha \right) =\frac{\widetilde{G}%
_{n+1,q}\left( \alpha \right) }{n+1}\text{.}
\end{equation*}

We now consider the function $\widetilde{G}_{q}\left( n:\alpha \right) $ as
the analytic continuation of weighted $q$-Genocchi numbers. All the weighted 
$q$-Genocchi numbers agree with $\widetilde{G}_{q}\left( n:\alpha \right) $,
the analytic continuation of weighted $q$-Genocchi numbers evaluated at $n$.
For $n\geq 0$, $\widetilde{G}_{q}\left( n:\alpha \right) =\widetilde{G}%
_{n,q}\left( \alpha \right) $.

We can now state $\widetilde{G}%
{\acute{}}%
_{q}\left( s:\alpha \right) $ in terms of $\widetilde{\zeta }%
{\acute{}}%
_{G,q}\left( s\mid \alpha \right) $, the derivative of $\widetilde{\zeta }%
_{G,q}\left( s:\alpha \right) $%
\begin{equation*}
\frac{\widetilde{G}_{q}\left( s+1:\alpha \right) }{s+1}=\widetilde{\zeta }%
_{G,q}\left( -s\mid \alpha \right) \text{, }\frac{\widetilde{G}%
{\acute{}}%
_{q}\left( s+1:\alpha \right) }{s+1}=\widetilde{\zeta }%
{\acute{}}%
_{G,q}\left( -s\mid \alpha \right) \text{.}
\end{equation*}

For $n,\alpha \in 
\mathbb{N}
\tbigcup \left\{ 0\right\} $ 
\begin{equation*}
\text{ }\frac{\widetilde{G}%
{\acute{}}%
_{q}\left( 2n+1:\alpha \right) }{2n+1}=\widetilde{\zeta }%
{\acute{}}%
_{G,q}\left( -2n\mid \alpha \right) \text{.}
\end{equation*}

This is suitable for the differential of the functional equation and so
supports the coherence of \noindent $\widetilde{G}_{q}\left( s:\alpha
\right) $ and $\widetilde{G}%
{\acute{}}%
_{q}\left( s:\alpha \right) $ with $\widetilde{G}_{n,q}\left( \alpha \right) 
$ and $\widetilde{\zeta }_{G,q}\left( s\mid \alpha \right) $. From the
analytic continuation of weighted $q$-Genocchi numbers, we derive as
follows: 
\begin{equation*}
\frac{\widetilde{G}_{q}\left( s+1:\alpha \right) }{s+1}=\widetilde{\zeta }%
_{G,q}\left( -s\mid \alpha \right) \text{ and }\frac{\widetilde{G}_{q}\left(
-s+1:\alpha \right) }{-s+1}=\widetilde{\zeta }_{G,q}\left( s\mid \alpha
\right) \text{.}
\end{equation*}

Moreover, we derive the following:

For $n\in 
\mathbb{N}
-\left\{ 1\right\} $%
\begin{equation*}
\frac{\widetilde{G}_{-n+1,q}\left( \alpha \right) }{-n+1}=\frac{\widetilde{G}%
_{q}\left( -n+1:\alpha \right) }{-n+1}=\widetilde{\zeta }_{G,q}\left( n\mid
\alpha \right) \text{.}
\end{equation*}

The curve $\widetilde{G}_{q}\left( s:a\right) $ review quickly the points $%
\widetilde{G}_{-s,q}\left( \alpha \right) $ and grows $\sim n$
asymptotically $\left( -n\right) \rightarrow -\infty $. The curve $%
\widetilde{G}_{q}\left( s:a\right) $ review quickly the point $\widetilde{G}%
_{q}\left( -s:a\right) $. Then, we procure the following:%
\begin{eqnarray*}
\lim_{n\rightarrow \infty }\frac{\widetilde{G}_{q}\left( -n+1:\alpha \right) 
}{-n+1} &=&\lim_{n\rightarrow \infty }\widetilde{\zeta }_{G,q}\left( n\mid
\alpha \right) =\lim_{n\rightarrow \infty }\left( \left[ 2:q\right]
\sum_{m=1}^{\infty }\frac{\left( -1\right) ^{m}q^{m}}{\left[ m:q^{\alpha }%
\right] ^{n}}\right) \\
&=&\lim_{n\rightarrow \infty }\left( -q\left[ 2:q\right] +\left[ 2:q\right]
\sum_{m=2}^{\infty }\frac{\left( -1\right) ^{m}q^{m}}{\left[ m:q^{\alpha }%
\right] ^{n}}\right) =-q^{2}\left[ 2:q^{-1}\right] \text{.}
\end{eqnarray*}%
From this, we easily note that

\begin{equation*}
\frac{\widetilde{G}_{q}\left( -n+1:\alpha \right) }{-n+1}=\widetilde{\zeta }%
_{G,q}\left( n\mid \alpha \right) \mapsto \frac{\widetilde{G}_{q}\left(
-s+1:\alpha \right) }{-s+1}=\widetilde{\zeta }_{G,q}\left( s\mid \alpha
\right) \text{.}
\end{equation*}

\section{\textbf{ANALYTIC CONTINUATION OF WEIGHTED }$q$\textbf{-GENOCCHI
POLYNOMIALS}}

For coherence with the redefinition of $\widetilde{G}_{n,q}\left( \alpha
\right) =\widetilde{G}_{q}\left( n:\alpha \right) $, we have%
\begin{equation*}
\widetilde{G}_{n,q}\left( x\mid \alpha \right) =q^{-\alpha x}\sum_{k=0}^{n}%
\binom{n}{k}q^{\alpha kx}\widetilde{G}_{k,q}\left( \alpha \right) \left[
x:q^{\alpha }\right] ^{n-k}\text{.}
\end{equation*}

Let $\Gamma \left( s\right) $ be Euler-gamma function. Then the analytic
continuation can be get as 
\begin{eqnarray*}
n &\mapsto &s\in 
\mathbb{R}
\text{, }x\mapsto w\in 
\mathbb{C}
\text{,} \\
\widetilde{G}_{n,q}\left( \alpha \right) &\mapsto &\widetilde{G}_{q}\left(
k+s-\left[ s\right] :\alpha \right) =\widetilde{\zeta }_{G,q}\left( -\left(
k+s-\left[ s\right] \right) \mid \alpha \right) \text{,} \\
\binom{n}{k} &=&\frac{\Gamma \left( n+1\right) }{\Gamma \left( n-k+1\right)
\Gamma \left( k+1\right) }\mapsto \frac{\Gamma \left( s+1\right) }{\Gamma
\left( 1+k+\left( s-\left[ s\right] \right) \right) \Gamma \left( 1+\left[ s%
\right] -k\right) } \\
\widetilde{G}_{s,q}\left( w\mid \alpha \right) &\mapsto &\widetilde{G}%
_{q}\left( s,w:\alpha \right) =q^{-\alpha w}\sum_{k=-1}^{\left[ s\right] }%
\frac{\Gamma \left( s+1\right) \widetilde{G}_{q}\left( k+\left( s-\left[ s%
\right] \right) :\alpha \right) q^{\alpha w\left( k+\left( s-\left[ s\right]
\right) \right) }}{\Gamma \left( 1+k+\left( s-\left[ s\right] \right)
\right) \Gamma \left( 1+\left[ s\right] -k\right) }\left[ w:q^{\alpha }%
\right] ^{\left[ s\right] -k} \\
&=&q^{-\alpha w}\sum_{k=0}^{\left[ s\right] +1}\frac{\Gamma \left(
s+1\right) \widetilde{G}_{q}\left( -1+k+\left( s-\left[ s\right] \right)
:\alpha \right) q^{\alpha w\left( k-1+\left( s-\left[ s\right] \right)
\right) }}{\Gamma \left( k+\left( s-\left[ s\right] \right) \right) \Gamma
\left( 2+\left[ s\right] -k\right) }\left[ w:q^{\alpha }\right] ^{\left[ s%
\right] +1-k}\text{.}
\end{eqnarray*}

Here $\left[ s\right] $ gives the integer part of s, and so $s-\left[ s%
\right] $ gives the fractional part.

Deformation of the curve $\widetilde{G}_{q}\left( 1,w:\alpha \right) $ into
the curve of $\widetilde{G}_{q}\left( 2,w:\alpha \right) $ is by means of
the real analytic cotinuation $\widetilde{G}_{q}\left( s,w:\alpha \right) $, 
$1\leq s\leq 2$, $-0.5\leq w\leq 0.5$.


\begin{thebibliography}{9}
\bibitem{Kim 1} T. Kim, Analytic continuation of $q$-Euler numbers and
polynomials, Applied Mathematics Letters 21 (2008) 1320-1323.

\bibitem{Kim 2} T. Kim, On explicit formulas of $p$-adic $q$-$L$-functions,
Kyushu J. Math. 43 (1994) 73--86.

\bibitem{kim 3} T. Kim, On $p$-adic interpolating function for $q$-Euler
numbers and its derivatives, J. Math. Anal. Appl. 339 (2008) 598--608.

\bibitem{kim 4} T. Kim, On the $q$-extension of Euler and Genocchi numbers,
J. Math. Anal. Appl. 326 (2007) 1458--1465.

\bibitem{Kim 5} T. Kim, On a $q$-analogue of the $p$-adic $\log $ gamma
functions and related integrals, Journal of Number Theory 76 (1999) 320-329.

\bibitem{Kim 6} T. Kim, On the analogs of Euler numbers and polynomials
associated with $p$-adic $q$-integral on $%
\mathbb{Z}
_{p}$ at $q=-1$, J. Math. Anal. Appl. 331 (2007) 779--792.

\bibitem{Araci 1} S. Araci, D. Erdal, J. J. Seo, A study on the fermionic $p$%
-adic $q$-integral representation on $%
\mathbb{Z}
_{p}$ associated with weighted $q$-Bernstein and $q$-Genocchi polynomials,
Abstract and Applied Analysis, Volume 2011, Article ID 649248, 10 pages.

\bibitem{Araci 2} S. Araci, M. Acikgoz and J. J. Seo, A study on the
weighted $q$-Genocchi numbers and polynomials with their interpolation
function, Honam Mathematical Journal (in press).
\end{thebibliography}
\end{document}